\documentstyle[draft%
,amscd%
,syntonly%
,amssymb%
,psfig%
]{amsart}

\newtheorem{thm}{Theorem}   

\newtheorem{lem}{Lemma}

\newtheorem{rem}{Remark}
\newtheorem{defn}{Definition}

\begin{document}



\title[bridge numbers]
{Additivity of Bridge Numbers of Knots}
\author{Jennifer Schultens}
\address{Department of Math and CS \\
   Emory University \\
   Atlanta, GA 30322} 
  
\email{jcs@@mathcs.emory.edu}
\thanks{Research partially supported by NSF grant DMS-9803826} 
\maketitle

\begin{abstract}
We provide a new proof of the following results of H. Schubert: If $K$
is a satellite knot with companion $J$ and pattern $(\hat V, L)$ with
index $k$, then the bridge numbers satisfy the following: $b(K) \geq k
\cdot (b(J))$.  In addition, if $K$ is a composite knot with summands
$J$ and $L$, then $b(K) = b(J) + b(L) - 1.$
\end{abstract}

\vspace{2 mm}

In ``\"Uber eine numerische Knoteninvariante'' \cite{S}, Horst
Schubert proved that for a satellite knot $K$ with companion $J$ and
pattern of index $k$, bridge numbers satisfy the inequality $b(K) \geq
k \cdot (b(J))$.  He also proved that for a composite knot $K$ with summands
$J$ and $L$, the bridge numbers satisfy $b(K) = b(J) + b(L) - 1$.  His
investigation was motivated by the question as to whether a knot can
have only finitely many companions.  Together with the fact that the
only bridge number one knot is the unknot, his result showed that the
answer to this question is yes.  

Schubert's main result may be recovered by a much shorter proof.  This
shorter proof grew out of an endeavour to recast the problem within
the framework of the thin position of a knot.  This framework turns
out to be far more refined than necessary.  The proof here does not
employ the notion of thin position.  It does, however, rely heavily on
the idea of rearranging the order in which critical points occur to
suit one's purpose, an idea fundamental to the notion of thin position
of knots and $3$-manifolds.  In this way it differs dramatically
from Schubert's proof.  It also differs from Schubert's in that it
relies on the consideration of Morse functions on $S^3$ whose level
sets are spheres (except for the maximum and minimum) and their
induced foliations.  This streamlines the terminology and the
complexity of the argument.  Schubert's proofs of the results reproven
here involve 25 pages containing 15 lemmas which involve a
consideration of up to three cases.

I wish to thank both Ray Lickorish and Marty Scharlemann for
independently suggesting that a more modern proof of this theorem
would be a welcome addition to the existing literature and for helpful
conversations.  I also wish to thank the Department of Pure
Mathematics and Mathematical Statistics at the University of Cambridge
for its hospitality.

\vspace{2 mm}

In the following $K$ will always be a knot in $S^3$ and $h:S^3
\rightarrow {\bf R}$ a Morse function with exactly two critical
points.  This last assumption guarantees that $h$ induces a foliation
of $S^3$ by spheres, along with one maximum that we denote by $\infty$
and one minimum that we denote by $-\infty$.

\begin{defn}
If the minima of $h_K$ occur below all maxima of $h_K$, then we say
that $K$ is in \underline{bridge position} with respect to $h$.  The
\underline{bridge number of K}, $b(K)$, is the minimal number of
maxima required for $h_K$.  (Note that this number is independent of
whether or not we require $K$ to be in bridge position.  Indeed, if
$h_K$ has $n$ maxima, then the maxima of $h_K$ can be raised, and the
minima of $h_K$ lowered, to obtain a copy of $K$ in bridge position with
$n$ maxima.)
\end{defn}

\begin{defn}
Let $J$ be a knot in $S^3$ and denote a small closed regular
neighborhood of $J$ by $\tilde V$.  Let $\hat V$ be an unknotted solid
torus in $S^3$ containing a knot $L$.  A map of $\hat V$ into $V$ maps
$L$ onto a knot $K$.  We call $K$ a \underline{satellite knot}, $J$ a
\underline{companion} of $K$, $V$ the \underline{companion torus} of
$K$ with respect to $J$ and $(\hat V, L)$ the \underline{pattern} (of
$K$ with respect to $J$).  The least number of times which a meridian
disk of $V$ intersects $L$ is called the \underline{index} of the
pattern.  (It is also called the wrapping number.)

In the special case in which the index of the pattern is 1, this
construction yields the connected sum of $J$ and $L$, and $V$ is
called a \underline{swallow-follow torus}.
\end{defn}

\begin{defn}
Suppose that $K$ is homotopically nontrivially contained in a
solid torus $V$.  Set $T = \partial V$.  Then $V$ is
\underline{taut with respect to $b(K)$}, if the number of
critical points of $h_T$ is minimal subject to the condition that
$h_K$ has $b(K)$ maxima.  
\end{defn}

\begin{defn}
Consider the singular foliation, $\cal F$$_T$, of $T$ induced by
$h_T$.  Let $\sigma$ be a leaf corresponding to a saddle singularity.
Then $\sigma$ consists of two circles, $s_1, s_2$, wedged at a point.
If either $s_1$ or $s_2$ is inessential in $T$, then we call $\sigma$
an \underline{inessential saddle}.  Otherwise, $\sigma$ is an
\underline{essential saddle}.
\end{defn}

\begin{lem} \label{popover}
(The Pop Over Lemma) Let $h, K, V, \cal F$$_T$ be as above. If $\cal
F$$_T$ contains inessential saddles, then, after an isotopy of $T$
that does not change $b(K)$ or the number of critical points of $h_T$,
there is an inessential saddle $\sigma$ in $\cal F$$_T$ for which the
following conditions hold:

\vspace{2 mm}

\noindent
1) $s_1$ bounds a disk $D_1 \subset T$ such that $\cal F$$_T$
restricted to $D_1$ contains only disjoint circles and one maximum or
minimum; and

\vspace{2 mm}

\noindent
2) for $L$ the level surface of $h$ containing $\sigma$, $D_1$
cobounds a $3$-ball $B$ with a disk $\tilde D_1 \subset L - T$, such
that $B$ does not contain $\infty$ or $-\infty$, and such that $s_2$
does not meet $B$ (i.e., such that $s_2$ lies outside of $\tilde
D_1$).
\end{lem}

\begin{pf}
The first condition on $\sigma$ may be satisfied by choosing $\sigma$
to be an inessential saddle in $\cal F$$_T$ that is innermost in $T$.
In this case $L - \partial D_1$ consists of two disks, $\hat D_1$ and
$\hat D_2$.  Together with $D_1$, both $\hat D_1$ and $\hat D_2$
cobound $3$-balls $\hat B_1, \hat B_2$, respectively.  One of these
$3$-balls, say $\hat B_2$, contains either $\infty$ or $-\infty$ and
the other contains neither.

If $s_2 \subset \hat D_2$, we may take $B = \hat B_1$, so suppose $s_2
\subset \hat D_1$.  Without loss of generality, we may assume that the
critical point of $D_1$ is a maximum.  In this case, consider a
monotone arc $\alpha$ disjoint from $K$, beginning at the maximum of
$D_1$, passing only through maxima of $T$ and ending at $\infty$.  Let
$a_1, \dots, a_n$ be the points at which $\alpha$ meets $T$, with
$a_n$ the highest such point.

Let $\beta$ be the subarc of $\alpha$ between $a_n$ and $\infty$ and
let $C'$ be a collar neighborhood of $\beta$.  After a small isotopy,
$T \cap C'$ consists of a small disk $D = a_n \times disk \subset T$.
Let $C''$ be a small $3$-ball centered at $\infty$ that is disjoint
from $T$.  Set $C = C' \cup C''$ and consider $T' = (T - D) \cup
(\partial C - D)$.  This describes an isotopy of $T$ that replaces
$\hat B_1$ by $\hat B_1 \cup C$ and replaces $\hat B_2$ by $\hat B_2 -
C$.  After a small tilt which turns $h_{T'}$ into a Morse function,
the maximum $a_n$ of $h_T$ has turned into a maximum of $h_{T'}$ at a
higher level.  No critical points need have been introduced for $h_K$
and the number of critical points of $h_{T'}$ is the same as that of
$h_T$.

By induction, we may assume that $\alpha$ is disjoint from $T$ except
at its initial point.  Then if $s_2 \subset \hat D_1$, this same
construction using $\beta = \alpha$ describes an isotopy of $T$
augmenting $\hat B_1$ to contain $\infty$ and shrinking $\hat B_2$ to
exclude $\infty$ without introducing any critical points of $h_K$ or
$h_T$.  We may then choose $B$ to be the shrunk version of $\hat B_2$.
\end{pf}

\begin{center}
fig. 1
\end{center}

\begin{lem} \label{popout}
(The Pop Out Lemma) Let $h, K, V, \cal F$$_T$ be as above.
If $V$ is taut with respect to $b(K)$, then there are no inessential
saddles in $\cal F$$_T$.
\end{lem}

\begin{pf}
Suppose there are inessential saddles.  Alter $T$ as in Lemma
\ref{popover} so that there is an inessential saddle $\sigma$
satisfying the conclusions of Lemma \ref{popover}.  We may assume that
$D_1$ contains a maximum and lies above $L$.  (The other case is
analogous.)  Here $(K \cup T) \cap int(B)$ can be shrunk horizontally
and lowered via an isotopy to lie just below $\tilde D_1$ (and above
any critical points of $h_K$ or $h_T$ below $\tilde D_1$).  This does
not change the nature or number of critical points of $h_K$ or $h_T$.

Now $D_1 \subset T$ can be replaced by $\tilde D_1$ to obtain $\tilde
T$.  After a small tilt, $\tilde T$ bounds a solid torus $\tilde V$
containing a copy of $K$ with $b(K)$ maxima, and $\tilde T$ is
isotopic to $T$, yet $h_{\tilde T}$ has two fewer critical points than
$h_T$.  (A maximum and an inessential saddle have been cancelled).
This contradicts the assumption that $V$ is taut with respect to
$b(K)$.
\end{pf}

\begin{center}
fig. 2
\end{center}

\begin{rem} \label{c3disk}
Consider a bicollar of an essential saddle $\sigma$ in $\cal F$$_T$.
It has three boundary components, $c_1, c_2, c_3$, where $c_i$ is
parallel to $s_i$ for $i = 1,2$.  Since $\chi(T) = 0$, it follows that
$c_3$ bounds a disk.  If there are no inessential saddles, then the
disk bounded by $c_3$ contains exactly one singular point, a maximum
or minimum.  We consider this maximum or minimum, $m_{\sigma}$, to be
the maximum or minimum corresponding to $\sigma$.

Conversely, if there are no inessential saddles in $\cal F$$_T$, then
every maximum or minimum corresponds to a saddle in this way, since
$\chi(T) = 0$.
\end{rem}

\begin{defn}
Let $\sigma, c_1, c_2, c_3$ be as above.  We may assume that $c_1$ and
$c_2$ are in the same level surface $L$ of $h$.  Then since $L$ is a
sphere, $c_1$ and $c_2$ cobound an annulus in $L$.  If a collar of
$c_1 \cup c_2$ in this annulus is contained in $V$, then $\sigma$ is a
\underline{nested saddle}.
\end{defn}

\begin{lem} \label{no nesting}
Let $h, K, V, \cal F$$_T$ be as above.  If $V$ is taut with respect to
$b(K)$, then $\cal F$$_T$ has no nested saddles.
\end{lem}

\begin{pf}
Suppose that there are nested saddles in $\cal F$$_T$.  

\vspace{1 mm}
Claim: There are also saddles in $\cal F$$_T$ that are not nested.

\vspace{1 mm} Let $\sigma$ be the highest saddle in $\cal F$$_T$.  For
$c_1,$ $c_2,$ $L$ as above, let $\hat D_1,$ $\hat D_2$ be the
(disjoint) disks bounded by $c_1,$ $c_2$ in $L$.  As $\sigma$ is the
highest saddle in $\cal F$$_T$, any curve in $T \cap interior(\hat
D_i)$ bounds a disk lying above $L$.  This implies that $\hat D_i$ is
isotopic to a disk whose interior is disjoint from $T$, i.e., lies
either entirely in $V$ or entirely in $S^3 - V$.  Since $c_i$ is
parallel to $s_i$, $c_i$ is essential in $T$.  Furthermore, since $V$
is knotted, $T$ is incompressible, whence $c_i$ is essential in the
closure of $S^3 - V$.  This implies that $\hat D_i$ must be isotopic
to a disk whose interior lies entirely in $V$ (in particular, $\hat
D_i$ is a meridian disk).  Thus $\sigma$ is not nested.

\vspace{1 mm}

If there are both saddles that are nested and saddles that are not
nested, then there must be an ``adjacent'' pair $\sigma_1, \sigma_2$
of essential saddles in $T$ with $\sigma_1$ nested, $\sigma_2$ not
nested, where ``adjacent'' means that one component, say $C$, of $T -
(\sigma_1 \cup \sigma_2)$ contains no critical points of $h_T$.
Consider the circles $s_1^i, s_2^i$ whose wedge is $\sigma_i$.
Without loss of generality, $s_1^1$ and $s_1^2$ meet $C$.

Again without loss of generality, we may assume that $\sigma_1$ lies
above $\sigma_2$ and hence that the component of $T - \sigma_1$ lying
above $\sigma_1$ and meeting both $s_1^1$ and $s_2^1$ is a disk
$D_3^1$. Construct a disk $D$ by adding $D_3^1$ to $C$ and capping off
$s_2^1$ with a level disk (a component of $h^{-1}(h(\sigma_1)) -
\sigma_1$).  Note that by the discussion above, this latter horizontal
portion of $D$ meets $K$ and $T$.

We now proceed as in Lemma \ref{popover} and Lemma \ref{popout}.  Here
$\partial D  = s_1^2$,  so $\partial D$  divides $h^{-1}(h(\sigma_2))$
into  two disks,  $\hat D_1$  and $\hat  D_2$, that  cobound $3$-balls
$\hat B_1$  and $\hat B_2$ together  with $D$.  By the  proof of Lemma
\ref{popover},  we may assume  that $\hat  B_2$ contains  $\infty$ and
that $s_2^2  \subset \hat B_2$.   We may thus shrink  horizontally and
lower $(K \cup T) \cap B$  as in the proof of Lemma \ref{popout}.  The
difference  is that here  $K \cup  T$ meets  $D$ along  its horizontal
portion.   As  $(K  \cup  T)  \cap  B$  is  shrunk  horizontally,  the
horizontal portion of $D$  is lowered while remaining horizontal.  The
portion   of   $B$  lying   above   $h^{-1}(h(\sigma_1))$  is   shrunk
horizontally as  necessary.  In the end, a  product neighborhood below
the original horizontal portion of $D$ ends up intersecting $K \cup T$
in vertical arcs and surfaces.

As in the proof of Lemma \ref{popout}, the number of critical points
of $h_T$ can be reduced without altering the number of critical points
of $h_K$, contradicting the fact that $V$ is taut with respect to
$b(K)$.
\end{pf}

\begin{center}
fig. 3
\end{center}

\begin{rem} \label{fat knot}
If $V$ is a knotted solid torus that is taut with respect to $b(K)$
then all saddles are essential and there are no nested saddles.  It
follows that if $L = h^{-1}(r)$ for some regular value $r$, then $V
\cap L$ consists of disks.  More specifically, let $\sigma_1,$
$\dots,$ $\sigma_n$ be the saddles in $\cal F$$_T$, and let $L_i =
h^{-1}(h(\sigma_i))$.  Recall that each saddle $\sigma$ corresponds to
a maximum or minimum $m_{\sigma}$ of $h_T$.  Between the level
surfaces $h^{-1}(h(\sigma))$ and $h^{-1}(h(m_{\sigma}))$ lies a
portion $B_{\sigma}$ of $V$ that is a $3$-ball.  Here $L_1 \cup \dots
\cup L_n$ cuts $V$ into $B_{\sigma_1},$ $\dots,$ $B_{\sigma_n}$ and
vertical cylinders.
\end{rem}

\begin{thm}
Suppose $K$ is a satellite knot with companion $J$, companion torus
$\tilde V$, pattern $(\hat V, L)$ and index $k$.  Then $b(K) \geq k
\cdot b(J)$.  In addition, if $K$ is the connected sum of two knots
$K_1$ and $K_2$, then $b(K) = b(K_1) + b(K_2) - 1$.
\end{thm}

\begin{pf}
We may assume that $V$ is taut with respect to $b(K)$.  Then $V$ is as
described in Remark \ref{fat knot}.  We obtain a Morse function on
$(S^3, J)$ by making $V$ very thin.  So $b(J)$ is less than or equal
to the number of maxima of $h_{T = \partial V}$.

Consider a maximum of $T$.  It corresponds to a saddle $\sigma$, where
$\sigma$ is the wedge of the circles $s_1, s_2$, bounding level
meridian disks $\tilde D_1, \tilde D_2$ of $V$.  Here $\tilde D_1 \cup
\tilde D_2$ cuts off a $3$-ball $B_{\sigma}$ as in Remark \ref{fat
knot}.  For distinct saddles $\sigma_i$ and $\sigma_j$, $B_{\sigma_i}$
and $B_{\sigma_j}$ are disjoint.  Since at least $k$ strands pass
through both $\tilde D_1$ and $\tilde D_2$, there are at least $k$
maxima of $K$ in $B_{\sigma}$.  Whence $b(K) \geq k \cdot b(J)$.

In the special case where $K$ is the connected sum of $K_1$ and $K_2$,
the sattelite construction may be used with $K_1$ the companion and
$(\hat V, K_2)$ the pattern (of index $1$).  By renumbering, if
necessary, we may assume that $b(K_1) \geq b(K_2)$.  Then we still
obtain a Morse function on $(S^3, K_1)$ as above.  Furthermore, if,
for each maximum of $T$ and $\sigma$, $B_{\sigma}$, $\tilde D_1$,
$\tilde D_2$ as above, $\vline \tilde D_i \cap K \vline \geq 2$ for $i
= 1, 2$, then $B_{\sigma}$ contains at least two maxima of $K$.
Hence $b(K) \geq 2 \cdot b(K_1) \geq b(K_1) + b(K_2) - 1$.  Thus we
may assume that there is a meridian disk contained in a level
surface of $h$ that intersects $K$ once.

Recall from Remark \ref{fat knot} that $V$ is comprised of
$B_{\sigma_1},$ $\dots,$ $B_{\sigma_n}$ and vertical cylinders.  In a
vertical cylinder, a meridian disk contained in a level surface that
intersects $K$ once may be used to move all critical points upwards or
downwards and out of the cylinder.  Thus the intersection of $K$ with
the cylinder becomes a monotone arc and the number of critical points
of $h_K$ is unchanged.

In $B_{\sigma_1}$, assume that $\vline \tilde D_1 \cap K \vline = 1$.
We may assume that $\sigma_1$ corresponds to a maximum of $h_T$.  Let
$\alpha$ be the subarc of $B_{\sigma_1} \cap K$ that connects $\tilde
D_1 \cap K$ to the closest maximum of $K$ in $B_{\sigma_1}$.  We may
assume that this maximum is the highest maximum of $K$ in
$B_{\sigma_1}$.  Then consider a disk $E$ in $B_{\sigma_1}$ for which
$\partial E$ consists of four subarcs: $\alpha,$ $a_1,$ $a_2,$ $a_3$,
where $a_1$ and $a_3$ are horizontal arcs connecting the endpoints of
$\alpha$ to $\partial B_{\sigma_1} \subset T$, and $a_2$ is an arc in
$\partial B_{\sigma_1}$ connecting the other endpoints of $a_1$ and
$a_3$, that runs over the maximum of $\partial B_{\sigma_1}$, and has
no other critical points.  We further require that $E \cap T = a_2$.

Claim: After an isotopy that does not change the number of critical
points of $h_K$, $E \cap K = \alpha$.  

Let $p_1, \dots, p_k$ be the points in $E \cap K - \alpha$ with $p_n$
the highest such point.  A small monotone subarc $\beta$ of $K$
containing $p_n$ may be replaced by a monotone arc $\beta'$ that
begins at one endpoint of $K - \beta$, travels parallel to $E$ until
it almost reaches $\partial B_{\sigma_1}$, then circles around to the
other side of $E$ along $\partial B_{\sigma_1}$ and travels parallel
to $E$ on the other side of $E$ until it meets the other endpoint of
$K - \beta$.  See fig. 4.  The result is isotopic to $K$ and has the
same number of critical points as $K$, yet one fewer intersection with
$E$.  The Claim follows by induction.

\begin{center}
fig. 4
\end{center}

Now $B_{\sigma_1} \cap K$ may be isotoped horizontally and downward,
so that after the isotopy this intersection consists of one arc with
exactly one critical point.  In the vertical solid cylinder meeting
$B_{\sigma_1}$ at $\tilde D_2$, $\vline \tilde D_2 \cap K \vline = 1$
and $\tilde D_2$ may be used to isotope $K$ so that all critical
points are moved to $B_{\sigma_2}$ (after relabelling) and $K$
intersects the solid cylinder in a single monotone arc.  Here $\sigma_2$
corresponds to a minimum, but an identical argument shows how proceed.
After a finite number of iterations of this procedure, $B_{\sigma_i}
\cap K$ consists of a single arc with one critical point for $i = 1,
\dots, n-1$ and $K$ intersects all cylindrical portions of $V$ in
monotone arcs.  Then, (since $\sigma_n$ corresponds to a minimum)
$\partial B_{\sigma_2}$ cuts $K$ into $K_1 - (subarc \; containing\;
a\; minimum)$ and $K_2 - (subarc \; containing\; a\; maximum)$.  This
proves that $b(K_1 \# K_2) \geq b(K_1) + b(K_2) - 1$.

The other inequality follows by considering a copy of $K_2$ in bridge
position realizing $b(K_2)$ lying below a copy of $K_1$ in bridge
position realizing $b(K_1)$ and taking the connected sum.
\end{pf}

\begin{figure}[hbtp]
\centering
\psfig{figure=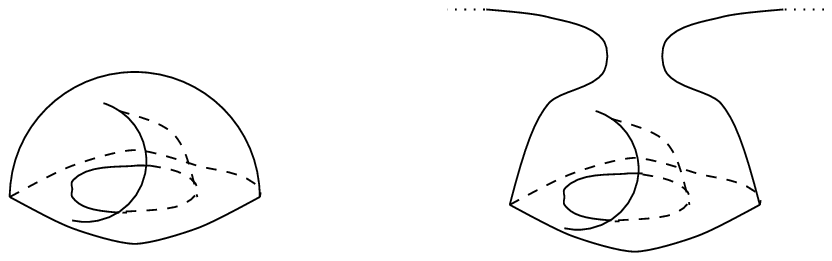}
\vspace{5 mm}
\begin{center}
fig. 1
\end{center}
\end{figure}
\vspace{1 mm}

\begin{figure}[hbtp]
\centering
\psfig{figure=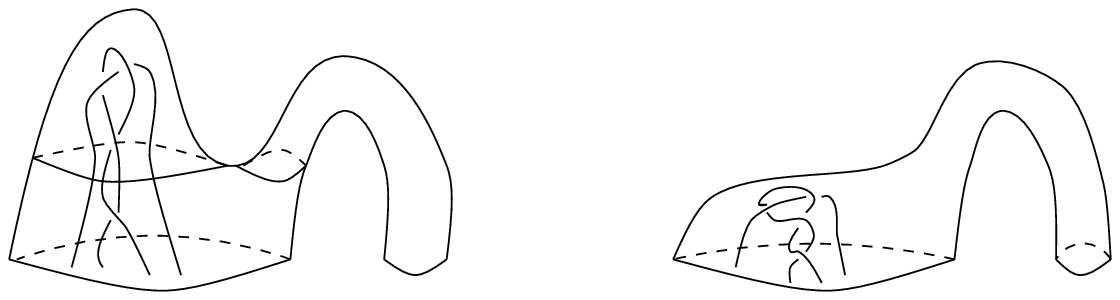}
\vspace{5 mm}
\begin{center}
fig. 2
\end{center}
\end{figure}
\vspace{1 mm}

\begin{figure}[hbtp]
\centering
\psfig{figure=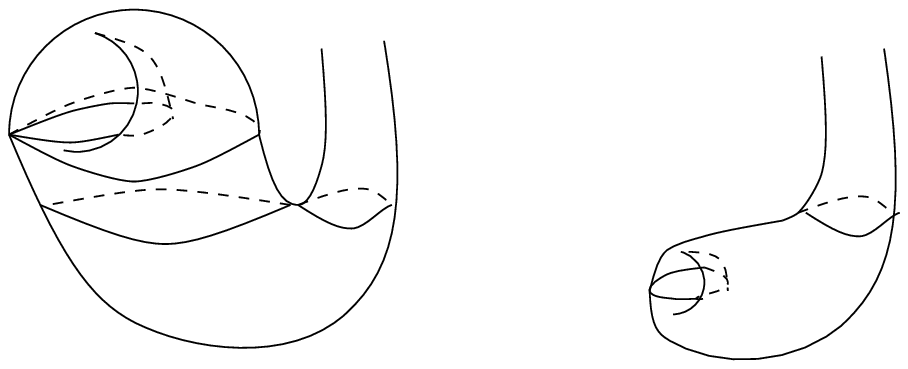}
\vspace{5 mm}
\begin{center}
fig. 3
\end{center}
\end{figure}
\vspace{1 mm}

\begin{figure}[hbtp]
\centering
\psfig{figure=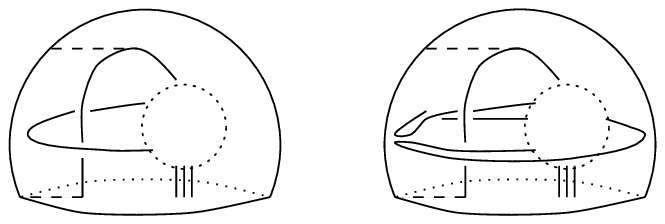}
\vspace{5 mm}
\begin{center}
fig. 4
\end{center}
\end{figure}
\vspace{1 mm}

\end{document}